\documentclass{amsproc}
\usepackage{color}

\newtheorem{theorem}{Theorem}[section]
\newtheorem{lemma}[theorem]{Lemma}
\newtheorem{proposition}[theorem]{Proposition}
\newtheorem{corollary}[theorem]{Corollary}
\theoremstyle{definition}

\newtheorem{assumption}[theorem]{Assumption}

\theoremstyle{remark}
\newtheorem{remark}[theorem]{Remark}

\newtheorem{notation}[theorem]{Notation}

\numberwithin{equation}{section}

\newcommand{\abs}[1]{\lvert#1\rvert}
\newcommand{\Abs}[1]{\Vert\lvert#1\rvert\Vert}
\newcommand{\blankbox}[2]{%
  \parbox{\columnwidth}{\centering

    \setlength{\fboxsep}{0pt}%
    \fbox{\raisebox{0pt}[#2]{\hspace{#1}}}%
  }%
}

\newcommand{\s}{{\rm supp}}
\newcommand{\R}{{\mathbb R}}
\newcommand{\RR}{\mathcal R}
\newcommand{\Z}{{\mathbb Z}}
\newcommand{\FF}{{\mathcal F}}
\newcommand{\GG}{{\mathcal G}}
\newcommand{\Int}{{{\rm Int}\,}}
\newcommand{\A}{{\mathbb A}}
\newcommand{\N}{{\mathbb N}}
\newcommand{\LL}{{\mathcal L}}
\newcommand{\IA}{{\stackrel{\circ}{\A}}}
\newcommand{\C}{{\mathbb C}}
\newcommand{\Fix}{{\rm Fix}}
\newcommand{\Moeb}{\mbox{M\"ob}_+(S^1_\C)}
\newcommand{\BB}{{\mathcal B}}
\newcommand{\PP}{{\mathcal P}}
\newcommand{\EE}{{\mathbb E}}
\newcommand{\HH}{{\mathcal H}}
\newcommand{\lev}{{\rm lev}}
\newcommand{\Stab}{{\rm Stab}}
\newcommand{\bd}{{\sc Proof}.\ }
\newcommand{\QQ}{{\mathcal Q}}
\newcommand{\VV}{{\mathbb V}}
\newcommand{\JJ}{{\mathcal J}}
\newcommand{\TT}{{\mathcal T}}
\newcommand{\II}{{\stackrel{\circ} I}}
\newcommand{\VVV}{{\mathcal V}}
\newcommand{\EEE}{{\mathcal E}}
\newcommand{\DD}{{\mathbb D}}
\newcommand{\nq}{\stackrel{\mbox{\Huge .}}{=}}

\begin{document}
\title[Kinematic expansive suspensions]
{Kinematic expansive suspensions of irrational rotations on the circle}


\author{Shigenori Matsumoto}
\address{Department of Mathematics, College of
Science and Technology, Nihon University, 1-8-14 Kanda, Surugadai,
Chiyoda-ku, Tokyo, 101-8308 Japan
}
\email{matsumo@math.cst.nihon-u.ac.jp
}
\thanks{The author is partially supported by Grant-in-Aid for
Scientific Research (C) No.\ 25400096.}


\date{\today }
\begin{abstract}
We shall show that the rotation of some irrational rotation number on
 the circle admits suspensions which are kinematic expansive.
\end{abstract}

\maketitle
\section{Introduction}

A continuous flow $\phi=\{\phi^t\}_{t\in\R}$ on a compact metric space
$X$ is called {\em kinematic expansive} if for any $\epsilon>0$, there
is $\delta(\epsilon)>0$ such that whenever $d(\phi^t(x),\phi^t(y))<\delta(\epsilon)$ for any
$t\in\R$, $y=\phi^s(x)$ for some $s\in(-\epsilon,\epsilon)$.

Given a homeomorphism $f$ of a compact metric space $Y$ and a continuous
fucntion $T:Y\to(0,\infty)$, we shall construct the  suspension flow
of $f$ with return time $T$ as follows. Let $\tilde\phi=\{\tilde\phi^t\}$
be the flow on $Y\times\R$ given by
$$
\tilde\phi^t(x,s)=(x,s+t).$$
Define a homeomorphism $F:Y\times\R\to Y\times\R$ by
$$
F(x,t)=(f(x),t-T(x)).$$
The infinite cyclic group $\langle F\rangle$  acts on $Y\times\R$
freely, properly discontinuously and cocompactly. The action commutes with
the flow $\tilde\phi$: 
$$
\tilde\phi^t\circ F=F\circ\tilde\phi^t,\ \ \forall t\in\R.
$$
Thus the flow $\tilde\phi$ induces a flow   on the quotient
space $\langle F\rangle\setminus (Y\times\R)$, which is called the
{\em suspension flow of $f$ with return time $T$}, 
denoted by ${\rm sus}(f,T)$. Notice that ${\rm
sus}(f,T)$
  admits a global cross section
$Y'$, the image of $Y\times\{0\}$ by the canonical
projection. The first return map of ${\rm sus}(f,T)$ with respect to $Y'$ is
$f$ and the return time is $T$. (Notice that a point $(x,0)$ on $Y'$
flows by time $T(x)$ to the point $(x,T(x))$, which is identified with
a point $(f(x),0)$ on $Y'$.) Of course the kinematic expansiveness of ${\rm sus}(f,T)$
strongly depends upon the choice of $T$.

In \cite{A}, A. Artigue studies among others  suspensions
of homeomorphisms of the circle $S^1$. He obtained:

\begin{theorem}
 Let $f$ be an orientation preserving nonminimal homeomorphism of
 $S^1$. Then
 $f$ admits a 
 kinematic expansive suspension if and only if there is a nonempty
 family $\{I_1,\cdots,I_r\}$ of finitely many nonempty open intervals such that
the wandering point set $W(f)$ of $f$ satisfies
$$
W(f)=\bigcup_{n\in\Z}\bigcup_{i=1}^rf^n(I_i).$$
\end{theorem} 

For the rotation
$R_\alpha$ by an irrational number $\alpha$, he showed that
if $T$ is absolutely continuous, then
${\rm sus}(R_\alpha,T)$ is not kinematic
expansive, and posed the problem for $T$ just continuous.   
The main result of the present paper is the following theorem.

\begin{theorem}\label{t2}
 There exist an irrational number $\alpha$ and a positive valued
 continuous function $T$ on $S^1$ such that ${\rm sus}(R_\alpha,T)$
is kinematic expansive.
\end{theorem}

\begin{remark}
 In fact, we have shown
for $\{\phi^t\}={\rm sus}(R_\alpha,T)$ in the above theorem that
for any $\epsilon>0$, there exists $\delta(\epsilon)>0$ such that whenever
$d(\phi^t(x),\phi^t(y))<\delta(\epsilon)$ for any {\em positive} $t$, then $y=\phi^s(x)$ for some
$s\in(-\epsilon,\epsilon)$. This is slightly stronger than the kinematic
 expansiveness.
\end{remark}
\bigskip
 {\sc Acknowledgement.}  Hearty thanks are due to
the referee for careful reading and helpful suggestions.

		       


\section{Proof of Theorem \ref{t2}}

We shall choose a particular irrational number
$\alpha\in(0,1)$ (explained later) and construct a return time map
$T:S^1\to(0,\infty)$ such that for some $\delta>0$,
\begin{equation}\label{e2}
y\neq x,\ \ \vert y-x\vert<\delta  \Longrightarrow\ \  \abs{T^{(n)}(y)-T^{(n)}(x)}>\delta,\ \ \exists n\in\N,
\end{equation}
where $\displaystyle T^{(n)}=\sum_{k=0}^{n-1}T\circ f^k$ is the $n$-th
return time.
This is sufficient for Theorem \ref{t2} since
we can choose the number
$\delta(\epsilon)$ in the definition of the kinematic expansiveness
as $\delta(\epsilon)=\min\{3^{-1}\epsilon,\delta\}$.
Actually we shall construct a {\em real valued} continuous function $T$
satisfying (\ref{e2}). We just need to add a positive
constant in order to make it positive valued.

For $x\in\R$, we denote its projected image on $S^1=\R/\Z$ by the same
letter $x$, and the distance in $S^1$ to $0$ by $\vert x\vert$. Notice
also that $S^1$ is an additive group so that for example $x-y$ makes
sense for $x,y\in S^1$.
We first prepare a useful lemma.

\begin{lemma}
 \label{l1}
Assume that for any $r\in(0,\delta)$, there is $x_r\in S^1$
 such that for some $m\in\N$,
\begin{equation*}
 \label{e3}
\abs{T^{(m)}(x_r+r)-T^{(m)}(x_r)}>3\delta.
\end{equation*} 
Then for any $x\in S^1$ and $y=x+r$,  the conclusion of (\ref{e2})
holds.
\end{lemma}

\bd
 For any $r\in(0,\delta)$, let $x_r$ and $m$ be as in the assumption of
Lemma \ref{l1}. For any $x\in S^1$, there exists $q\in\N$ such that $R_\alpha^q(x)$ is
arbitrarily close to $x_r$ (and thus $R_\alpha^q(x+r)$ to $x_r+r$). By the uniform continuity of $T^{(m)}$, one may assume
\begin{equation*}
 \label{e4}
\abs{T^{(m)}(R_\alpha^q(x))-T^{(m)}(x_r)}<\delta/2,\ \ 
\abs{T^{(m)}(R_\alpha^q(x+r))-T^{(m)}(x_r+r)}<\delta/2.
\end{equation*}
If $\abs{T^{(q)}(x+r)-T^{(q)}(x)}>\delta$, then there is nothing to prove.
Otherwise we have
$$
\abs{T^{(m+q)}(x+r)-T^{(m+q)}(x)}>\delta,$$
as is required. \qed

\bigskip
We choose the irrational number $\alpha$ by the continued fraction as
$$
\alpha=\frac{1}{a+\frac{1}{a+\frac{1}{a+\cdots}}},$$
where $a$ is an integer $\geq10^{10}$.
That is, $$\alpha=\frac{1}{2}(-a+\sqrt{a^2+4}).$$

In fact, the arguments in what follows work for much smaller value of
$a$. On the other hand, they are not applicable to Liouville numbers.
So we make the assumption $a\geq10^{10}$ in order to make
various estimates easier.

Let $p_n/q_n$ be the $n$-th convergent of $\alpha$. The denominator
$q_n$ is obtained inductively as:
\begin{equation*}
 \label{e5}
q_{n+2}=aq_{n+1}+q_{n},\ \ q_0=1,\ \ q_1=a.
\end{equation*}

Let $c$ be the positive solution of $x^2=ax+1$:
$$
c=2^{-1}(a+\sqrt{a^2+4}).$$
Thus $c$ is a number slightly
bigger than $a$. We have
$$q_n=Ac^n+B(-c^{-1})^n,\ \mbox{ where } \ 
A=\frac{a+\sqrt{a^2+4}}{2\sqrt{a^2+4}},\ \
B=\frac{-a+\sqrt{a^2+4}}{2\sqrt{a^2+4}}.$$
Thus $A$ and $B$ are positive numbers satisfying $A+B=1$, and $A$ is
almost 1. 
Since $c>a\geq 10^{10}$, we have
\begin{equation}
 \label{e6}
q_n\nq Ac^n.
\end{equation}

\begin{notation}
For $a,b>0$, $a\nq b$ means $a/b\in(\frac{100}{101},\frac{101}{100})$.
\end{notation}

\medskip
 It is well known that $q_n$ is the {\em closest return time} for the rotation
$R_\alpha:S^1\to S^1$. That is, $R_\alpha^{q_n}(x)$ is the closest to
$x$ among the points $R_\alpha(x),\cdots,R_\alpha^{q_n-1}(x),R_\alpha^{q_n}(x)$.
More precisely (letting $x=0$),
\begin{equation}
 \label{e7}
\vert q_n\alpha\vert<\vert i\alpha\vert,\ \ \forall
i\in\{1,2,\cdots,q_n-1\}.
\end{equation}
The point $q_n\alpha$ is very close to $0$, lies on the right side of 0 if
$n$ is odd, and on the left if $n$ is even. Let $I_n$ be the
smaller closed interval in $S^1$ bounded by $0$ and $q_n\alpha$. Consider the
first return map of $R_\alpha$ for the interval $I_n\cup I_{n+1}$.
The part $I_{n+1}$ returns to $I_n\cup I_{n+1}$ for the first time
by the $q_n$ iterate of $R_\alpha$, and the part $I_n\setminus\{0\}$ 
by the $q_{n+1}$ iterate. (All this follows from (\ref{e7}).)
See Figure 1 for even $n$.

\begin{figure}\caption{}
{\unitlength 0.1in%
\begin{picture}( 40.0000, 20.3000)( 22.9000,-29.8000)%
%
\special{pn 8}%
\special{pa 2370 1260}%
\special{pa 6290 1260}%
\special{fp}%
%
\special{pn 8}%
\special{pa 2370 2980}%
\special{pa 6250 2980}%
\special{fp}%
%
\special{pn 4}%
\special{sh 1}%
\special{ar 5170 1260 16 16 0  6.28318530717959E+0000}%
\special{sh 1}%
\special{ar 5170 1260 16 16 0  6.28318530717959E+0000}%
%
\special{pn 4}%
\special{sh 1}%
\special{ar 6270 1260 16 16 0  6.28318530717959E+0000}%
\special{sh 1}%
\special{ar 6270 1260 16 16 0  6.28318530717959E+0000}%
%
\special{pn 4}%
\special{sh 1}%
\special{ar 2380 1260 16 16 0  6.28318530717959E+0000}%
\special{sh 1}%
\special{ar 2380 1260 16 16 0  6.28318530717959E+0000}%
%
\special{pn 4}%
\special{sh 1}%
\special{ar 2360 2980 16 16 0  6.28318530717959E+0000}%
\special{sh 1}%
\special{ar 2350 2980 16 16 0  6.28318530717959E+0000}%
%
\special{pn 4}%
\special{sh 1}%
\special{ar 6260 2980 16 16 0  6.28318530717959E+0000}%
\special{sh 1}%
\special{ar 6260 2980 16 16 0  6.28318530717959E+0000}%
%
\special{pn 4}%
\special{pa 5240 1420}%
\special{pa 6180 2880}%
\special{fp}%
\special{sh 1}%
\special{pa 6180 2880}%
\special{pa 6161 2813}%
\special{pa 6151 2835}%
\special{pa 6127 2835}%
\special{pa 6180 2880}%
\special{fp}%
%
\special{pn 4}%
\special{pa 2490 1440}%
\special{pa 3330 2860}%
\special{fp}%
\special{sh 1}%
\special{pa 3330 2860}%
\special{pa 3313 2792}%
\special{pa 3303 2814}%
\special{pa 3279 2813}%
\special{pa 3330 2860}%
\special{fp}%
%
\special{pn 4}%
\special{pa 5100 1400}%
\special{pa 2510 2840}%
\special{dt 0.045}%
\special{sh 1}%
\special{pa 2510 2840}%
\special{pa 2578 2825}%
\special{pa 2557 2814}%
\special{pa 2559 2790}%
\special{pa 2510 2840}%
\special{fp}%
%
\special{pn 4}%
\special{pa 6200 1420}%
\special{pa 3540 2870}%
\special{dt 0.045}%
\special{sh 1}%
\special{pa 3540 2870}%
\special{pa 3608 2856}%
\special{pa 3587 2844}%
\special{pa 3589 2821}%
\special{pa 3540 2870}%
\special{fp}%
%
\special{pn 4}%
\special{sh 1}%
\special{ar 3450 2980 16 16 0  6.28318530717959E+0000}%
\special{sh 1}%
\special{ar 3450 2970 16 16 0  6.28318530717959E+0000}%
%
\put(52.0000,-10.2000){\makebox(0,0)[lb]{}}%
\put(51.1000,-10.9000){\makebox(0,0)[lb]{$0$}}%
\put(61.5000,-10.9000){\makebox(0,0)[lb]{$q_{n+1}\alpha$}}%
\put(22.9000,-10.9000){\makebox(0,0)[lb]{$q_n\alpha$}}%
\put(59.9000,-24.2000){\makebox(0,0)[lb]{$R_\alpha^{q_{n+1}}$}}%
\put(40.1000,-21.9000){\makebox(0,0)[lb]{$R_\alpha^{q_n}$}}%
\put(35.4000,-10.8000){\makebox(0,0)[lb]{$I_n$}}%
\put(55.8000,-10.9000){\makebox(0,0)[lb]{$I_{n+1}$}}%
\end{picture}}%

\end{figure}

The intervals $I_n,R_\alpha(I_n),\ldots,R_\alpha^{q_{n+1}-1}(I_n)$ and
$I_{n+1},R_\alpha(I_{n+1})$, $\ldots,R_\alpha^{q_n-1}(I_{n+1})$ yield 
a partition of $S^1$ (a covering of $S^1$ by nonoverlapping
intervals.) In fact, (\ref{e7})  implies that the intervals are
nonoverlapping. On the other hand, their total length is 1 by the
following well-known equality:
$$
q_{n+1}\vert I_n\vert+q_n\vert I_{n+1}\vert=\vert q_{n+1}(q_n\alpha-p_n)+q_n(p_{n+1}-q_{n+1}\alpha)\vert=\vert-q_{n+1}p_n+q_np_{n+1}\vert=1.$$

 Dynamically they form a Rochlin tower as is
depicted in Figure 2.

\begin{figure}\caption{}
{\unitlength 0.1in%
\begin{picture}( 42.2000, 34.0000)( 15.7000,-40.7000)%
%
\special{pn 8}%
\special{pa 1740 1040}%
\special{pa 4430 1040}%
\special{fp}%
\special{pa 4460 1040}%
\special{pa 4460 1040}%
\special{fp}%
%
\special{pn 8}%
\special{pa 1730 1330}%
\special{pa 4410 1330}%
\special{fp}%
%
\special{pn 8}%
\special{pa 1750 2460}%
\special{pa 5740 2460}%
\special{fp}%
%
\special{pn 8}%
\special{pa 1750 2760}%
\special{pa 5750 2760}%
\special{fp}%
%
\special{pn 8}%
\special{pa 1730 2140}%
\special{pa 4420 2140}%
\special{fp}%
%
\special{pn 8}%
\special{pa 1730 3520}%
\special{pa 5780 3520}%
\special{fp}%
%
\special{pn 8}%
\special{pa 1740 3860}%
\special{pa 5790 3860}%
\special{fp}%
%
\special{pn 8}%
\special{pa 1750 1050}%
\special{pa 1750 3870}%
\special{dt 0.045}%
%
\special{pn 8}%
\special{pa 4420 1040}%
\special{pa 4430 3850}%
\special{dt 0.045}%
%
\special{pn 8}%
\special{pa 5770 2480}%
\special{pa 5770 3850}%
\special{dt 0.045}%
%
\special{pn 8}%
\special{pa 2900 3780}%
\special{pa 2900 3870}%
\special{fp}%
%
\special{pn 4}%
\special{pa 3200 3020}%
\special{pa 3200 2900}%
\special{fp}%
\special{sh 1}%
\special{pa 3200 2900}%
\special{pa 3180 2967}%
\special{pa 3200 2953}%
\special{pa 3220 2967}%
\special{pa 3200 2900}%
\special{fp}%
%
\special{pn 4}%
\special{pa 3200 2090}%
\special{pa 3200 1920}%
\special{fp}%
\special{sh 1}%
\special{pa 3200 1920}%
\special{pa 3180 1987}%
\special{pa 3200 1973}%
\special{pa 3220 1987}%
\special{pa 3200 1920}%
\special{fp}%
%
\special{pn 4}%
\special{pa 3210 1590}%
\special{pa 3210 1420}%
\special{fp}%
\special{sh 1}%
\special{pa 3210 1420}%
\special{pa 3190 1487}%
\special{pa 3210 1473}%
\special{pa 3230 1487}%
\special{pa 3210 1420}%
\special{fp}%
%
\special{pn 4}%
\special{pa 3210 1280}%
\special{pa 3210 1090}%
\special{fp}%
\special{sh 1}%
\special{pa 3210 1090}%
\special{pa 3190 1157}%
\special{pa 3210 1143}%
\special{pa 3230 1157}%
\special{pa 3210 1090}%
\special{fp}%
%
\special{pn 4}%
\special{pa 3200 2400}%
\special{pa 3200 2230}%
\special{fp}%
\special{sh 1}%
\special{pa 3200 2230}%
\special{pa 3180 2297}%
\special{pa 3200 2283}%
\special{pa 3220 2297}%
\special{pa 3200 2230}%
\special{fp}%
%
\special{pn 4}%
\special{pa 3210 2680}%
\special{pa 3210 2530}%
\special{fp}%
\special{sh 1}%
\special{pa 3210 2530}%
\special{pa 3190 2597}%
\special{pa 3210 2583}%
\special{pa 3230 2597}%
\special{pa 3210 2530}%
\special{fp}%
%
\special{pn 4}%
\special{pa 3190 3450}%
\special{pa 3190 3280}%
\special{fp}%
\special{sh 1}%
\special{pa 3190 3280}%
\special{pa 3170 3347}%
\special{pa 3190 3333}%
\special{pa 3210 3347}%
\special{pa 3190 3280}%
\special{fp}%
%
\special{pn 4}%
\special{pa 3190 3800}%
\special{pa 3190 3600}%
\special{fp}%
\special{sh 1}%
\special{pa 3190 3600}%
\special{pa 3170 3667}%
\special{pa 3190 3653}%
\special{pa 3210 3667}%
\special{pa 3190 3600}%
\special{fp}%
%
\special{pn 4}%
\special{pa 5070 3810}%
\special{pa 5070 3620}%
\special{fp}%
\special{sh 1}%
\special{pa 5070 3620}%
\special{pa 5050 3687}%
\special{pa 5070 3673}%
\special{pa 5090 3687}%
\special{pa 5070 3620}%
\special{fp}%
%
\special{pn 4}%
\special{pa 5070 3450}%
\special{pa 5070 3300}%
\special{fp}%
\special{sh 1}%
\special{pa 5070 3300}%
\special{pa 5050 3367}%
\special{pa 5070 3353}%
\special{pa 5090 3367}%
\special{pa 5070 3300}%
\special{fp}%
%
\special{pn 4}%
\special{pa 5070 3040}%
\special{pa 5070 2820}%
\special{fp}%
\special{sh 1}%
\special{pa 5070 2820}%
\special{pa 5050 2887}%
\special{pa 5070 2873}%
\special{pa 5090 2887}%
\special{pa 5070 2820}%
\special{fp}%
%
\special{pn 4}%
\special{pa 5070 2700}%
\special{pa 5070 2550}%
\special{fp}%
\special{sh 1}%
\special{pa 5070 2550}%
\special{pa 5050 2617}%
\special{pa 5070 2603}%
\special{pa 5090 2617}%
\special{pa 5070 2550}%
\special{fp}%
%
\special{pn 4}%
\special{pa 4400 2530}%
\special{pa 1830 3800}%
\special{fp}%
\special{sh 1}%
\special{pa 1830 3800}%
\special{pa 1899 3788}%
\special{pa 1878 3776}%
\special{pa 1881 3753}%
\special{pa 1830 3800}%
\special{fp}%
%
\special{pn 4}%
\special{pa 5690 2550}%
\special{pa 3020 3810}%
\special{fp}%
\special{sh 1}%
\special{pa 3020 3810}%
\special{pa 3089 3800}%
\special{pa 3068 3787}%
\special{pa 3072 3763}%
\special{pa 3020 3810}%
\special{fp}%
%
\special{pn 4}%
\special{pa 1810 1140}%
\special{pa 2810 3800}%
\special{fp}%
\special{sh 1}%
\special{pa 2810 3800}%
\special{pa 2805 3731}%
\special{pa 2791 3750}%
\special{pa 2768 3745}%
\special{pa 2810 3800}%
\special{fp}%
%
\special{pn 4}%
\special{pa 4480 1140}%
\special{pa 5710 3750}%
\special{fp}%
\special{sh 1}%
\special{pa 5710 3750}%
\special{pa 5700 3681}%
\special{pa 5687 3702}%
\special{pa 5663 3698}%
\special{pa 5710 3750}%
\special{fp}%
\put(43.3000,-40.9000){\makebox(0,0)[lb]{$0$}}%
\put(56.0000,-40.6000){\makebox(0,0)[lb]{$q_{n+1}\alpha$}}%
\put(16.3000,-40.8000){\makebox(0,0)[lb]{$q_n\alpha$}}%
\put(50.7000,-42.8000){\makebox(0,0)[lb]{$I_{n+1}$}}%
\put(29.1000,-43.0000){\makebox(0,0)[lb]{$I_n$}}%
\put(27.2000,-40.7000){\makebox(0,0)[lb]{$(q_n+q_{n+1})\alpha$}}%
\put(45.1000,-24.1000){\makebox(0,0)[lb]{$(q_n-1)\alpha$}}%
\put(57.5000,-24.0000){\makebox(0,0)[lb]{$(q_n+q_{n+1}-1)\alpha$}}%
\put(15.7000,-9.4000){\makebox(0,0)[lb]{$(q_n+q_{n+1}-1)\alpha$}}%
\put(44.3000,-9.0000){\makebox(0,0)[lb]{$(q_{n+1}-1)\alpha$}}%
\put(50.7000,-20.6000){\makebox(0,0)[lb]{$R_\alpha^{q_n-1}(I_{n+1})$}}%
\put(51.5000,-34.8000){\makebox(0,0)[lb]{$R_\alpha(I_{n+1})$}}%
\put(33.0000,-34.6000){\makebox(0,0)[lb]{$R_\alpha(I_n)$}}%
\put(29.9000,-9.5000){\makebox(0,0)[lb]{$R_\alpha^{q_{n+1}-1}(I_n)$}}%
\end{picture}}%

\end{figure}

\bigskip Let us begin the construction of the function $T$. As
we remarked before, $T$ is to be positive and negative valued.
We shall construct $T$ as 
$T=\sum_{n=1}^\infty T_n$, where $T_n$ is a continuous function such
that
\begin{equation*}
 \label{e8}
\sum_{n=1}^\infty\vert T_n\vert_\infty<\infty.
\end{equation*}
To describe $T_n$, we shall first define a function $\chi_I:S^1\to\R$ for a 
interval $I=[r,s]$ of $S^1$ as follows. See Figure 3.
The Lipschitz constant of $\chi_I$ is 1.
$$
\chi_I(x)=\left\{ \begin{array}{ll}
x-r & \mbox{ if }r\leq x \leq 3^{-1}(2r+s) \\
3^{-1}(s-r) &  \mbox{ if }3^{-1}(2r+s)\leq x\leq 3^{-1}(r+2s) \\
-x+s & \mbox{ if }3^{-1}(r+2s)\leq x\leq s \\
 0  & \mbox{ for other $x$ }.\end{array}\right.$$

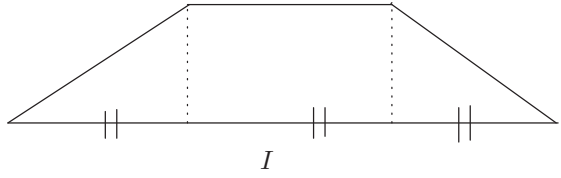
\begin{figure}\caption{The graph of $\chi_I$. The slope is 0, $\pm 1$.}
{\unitlength 0.1in%
\begin{picture}( 28.8000,  8.4000)( 23.3000,-18.3000)%
%
\special{pn 8}%
\special{pa 2330 1720}%
\special{pa 5210 1720}%
\special{fp}%
\put(36.5000,-19.6000){\makebox(0,0)[lb]{$I$}}%
%
\special{pn 8}%
\special{pa 3270 1720}%
\special{pa 3270 1090}%
\special{dt 0.045}%
%
\special{pn 8}%
\special{pa 4340 1720}%
\special{pa 4340 1090}%
\special{dt 0.045}%
\special{pa 4340 1090}%
\special{pa 4330 1100}%
\special{dt 0.045}%
%
\special{pn 8}%
\special{pa 3270 1100}%
\special{pa 4340 1100}%
\special{fp}%
%
\special{pn 8}%
\special{pa 2330 1720}%
\special{pa 3280 1100}%
\special{fp}%
%
\special{pn 8}%
\special{pa 4340 1100}%
\special{pa 5200 1720}%
\special{fp}%
%
\special{pn 8}%
\special{pa 2840 1660}%
\special{pa 2840 1800}%
\special{fp}%
%
\special{pn 8}%
\special{pa 2900 1650}%
\special{pa 2900 1800}%
\special{fp}%
%
\special{pn 8}%
\special{pa 3930 1640}%
\special{pa 3930 1800}%
\special{fp}%
%
\special{pn 8}%
\special{pa 3990 1640}%
\special{pa 3990 1790}%
\special{fp}%
%
\special{pn 8}%
\special{pa 4690 1640}%
\special{pa 4690 1820}%
\special{fp}%
%
\special{pn 8}%
\special{pa 4750 1630}%
\special{pa 4750 1810}%
\special{fp}%
\end{picture}}%

\end{figure}

For each $n$, define $j_n=\lfloor3^{-1}(q_{n}-1)\rfloor$ and
$j_n'=q_n-2j_n$, where $\lfloor x\rfloor$ denotes the integer part
of a real number $x$. For $n$ even, define
$$
T_n=\sum_{i=0}^{j_{n+1}-1}\chi_{R_\alpha^i(I_n)}-2^{-1}
\sum_{i=j_{n+1}'}^{q_{n+1}-1}\chi_{R_\alpha^i(I_n)}$$
Notice that $j_{n+1}'$ is slightly bigger than $j_{n+1}$:
in fact, $0\leq j'_{n+1}-j_{n+1}<4$.
The first term of $T_n$ has $j_{n+1}$ summands, while the second $2j_{n+1}$.
See Figure 4.

\begin{figure}\caption{The function $T_n$. The values at dotted points
 sum up to 0.}
{\unitlength 0.1in%
\begin{picture}( 28.9000, 54.8000)( 23.2000,-62.9000)%
%
\special{pn 8}%
\special{pa 2320 1180}%
\special{pa 4280 1180}%
\special{fp}%
%
\special{pn 8}%
\special{pa 2320 1180}%
\special{pa 3070 1440}%
\special{fp}%
%
\special{pn 8}%
\special{pa 2340 1660}%
\special{pa 4320 1660}%
\special{fp}%
%
\special{pn 8}%
\special{pa 2350 1660}%
\special{pa 3050 1880}%
\special{fp}%
%
\special{pn 8}%
\special{pa 3050 1880}%
\special{pa 3730 1880}%
\special{fp}%
%
\special{pn 8}%
\special{pa 3740 1880}%
\special{pa 4290 1660}%
\special{fp}%
%
\special{pn 8}%
\special{pa 3080 1440}%
\special{pa 3780 1440}%
\special{fp}%
%
\special{pn 8}%
\special{pa 3790 1440}%
\special{pa 4290 1180}%
\special{fp}%
%
\special{pn 8}%
\special{pa 2340 2650}%
\special{pa 4290 2650}%
\special{fp}%
%
\special{pn 8}%
\special{pa 3010 2930}%
\special{pa 3760 2930}%
\special{fp}%
%
\special{pn 8}%
\special{pa 3760 2940}%
\special{pa 4270 2650}%
\special{fp}%
%
\special{pn 8}%
\special{pa 2400 3140}%
\special{pa 4280 3140}%
\special{fp}%
%
\special{pn 8}%
\special{pa 2430 3140}%
\special{pa 3030 3390}%
\special{fp}%
%
\special{pn 8}%
\special{pa 3020 3390}%
\special{pa 3690 3390}%
\special{fp}%
%
\special{pn 8}%
\special{pa 3690 3400}%
\special{pa 4260 3160}%
\special{fp}%
%
\special{pn 8}%
\special{pa 2350 3600}%
\special{pa 4290 3600}%
\special{fp}%
%
\special{pn 8}%
\special{pa 3020 3840}%
\special{pa 3790 3840}%
\special{fp}%
%
\special{pn 8}%
\special{pa 3800 3840}%
\special{pa 4280 3600}%
\special{fp}%
%
\special{pn 8}%
\special{pa 2380 4630}%
\special{pa 4310 4630}%
\special{fp}%
%
\special{pn 8}%
\special{pa 2370 4630}%
\special{pa 2970 4120}%
\special{fp}%
%
\special{pn 8}%
\special{pa 2970 4140}%
\special{pa 3790 4140}%
\special{fp}%
%
\special{pn 8}%
\special{pa 3790 4150}%
\special{pa 4290 4620}%
\special{fp}%
%
\special{pn 8}%
\special{pa 2390 5460}%
\special{pa 5200 5460}%
\special{fp}%
%
\special{pn 8}%
\special{pa 2970 4950}%
\special{pa 3740 4950}%
\special{fp}%
%
\special{pn 8}%
\special{pa 3730 4960}%
\special{pa 4330 5450}%
\special{fp}%
%
\special{pn 8}%
\special{pa 2410 6300}%
\special{pa 5210 6300}%
\special{fp}%
%
\special{pn 8}%
\special{pa 2390 6300}%
\special{pa 2940 5810}%
\special{fp}%
%
\special{pn 8}%
\special{pa 2940 5820}%
\special{pa 3720 5820}%
\special{fp}%
%
\special{pn 8}%
\special{pa 3710 5830}%
\special{pa 4360 6300}%
\special{fp}%
%
\special{pn 4}%
\special{pa 3370 4010}%
\special{pa 4690 4010}%
\special{dt 0.045}%
\put(30.1000,-65.2000){\makebox(0,0)[lb]{$I_n$}}%
\put(46.1000,-65.0000){\makebox(0,0)[lb]{$I_{n+1}$}}%
\put(30.1000,-56.5000){\makebox(0,0)[lb]{$R_\alpha(I_n)$}}%
\put(45.5000,-56.3000){\makebox(0,0)[lb]{$R_\alpha(I_{n+1})$}}%
\put(29.7000,-10.4000){\makebox(0,0)[lb]{$R_\alpha^{q_{n+1}-1}(I_n)$}}%
\put(48.4000,-41.1000){\makebox(0,0)[lb]{around one third}}%
%
\special{pn 8}%
\special{pa 3380 3060}%
\special{pa 4730 3060}%
\special{dt 0.045}%
\put(48.4000,-31.5000){\makebox(0,0)[lb]{around one half}}%
%
\special{pn 8}%
\special{pa 2350 2660}%
\special{pa 3030 2930}%
\special{fp}%
%
\special{pn 8}%
\special{pa 2360 3590}%
\special{pa 3010 3840}%
\special{fp}%
%
\special{pn 8}%
\special{pa 2400 5460}%
\special{pa 2990 4950}%
\special{fp}%
%
\special{pn 4}%
\special{sh 1}%
\special{ar 2700 1180 16 16 0  6.28318530717959E+0000}%
\special{sh 1}%
\special{ar 2700 1180 16 16 0  6.28318530717959E+0000}%
%
\special{pn 4}%
\special{sh 1}%
\special{ar 2710 1660 16 16 0  6.28318530717959E+0000}%
\special{sh 1}%
\special{ar 2720 1660 16 16 0  6.28318530717959E+0000}%
%
\special{pn 4}%
\special{sh 1}%
\special{ar 2720 2660 8 8 0  6.28318530717959E+0000}%
\special{sh 1}%
\special{ar 2720 2660 8 8 0  6.28318530717959E+0000}%
%
\special{pn 4}%
\special{sh 1}%
\special{ar 2730 2650 16 16 0  6.28318530717959E+0000}%
\special{sh 1}%
\special{ar 2740 2650 16 16 0  6.28318530717959E+0000}%
%
\special{pn 4}%
\special{sh 1}%
\special{ar 2740 3140 16 16 0  6.28318530717959E+0000}%
\special{sh 1}%
\special{ar 2750 3140 16 16 0  6.28318530717959E+0000}%
%
\special{pn 4}%
\special{sh 1}%
\special{ar 2740 3610 16 16 0  6.28318530717959E+0000}%
\special{sh 1}%
\special{ar 2750 3590 16 16 0  6.28318530717959E+0000}%
%
\special{pn 4}%
\special{sh 1}%
\special{ar 2750 4630 16 16 0  6.28318530717959E+0000}%
\special{sh 1}%
\special{ar 2730 4630 16 16 0  6.28318530717959E+0000}%
%
\special{pn 4}%
\special{sh 1}%
\special{ar 2740 5460 16 16 0  6.28318530717959E+0000}%
\special{sh 1}%
\special{ar 2760 5460 16 16 0  6.28318530717959E+0000}%
%
\special{pn 4}%
\special{sh 1}%
\special{ar 2760 6300 16 16 0  6.28318530717959E+0000}%
\special{sh 1}%
\special{ar 2750 6300 16 16 0  6.28318530717959E+0000}%
\put(26.8000,-65.0000){\makebox(0,0)[lb]{$x$}}%
\end{picture}}%

\end{figure}
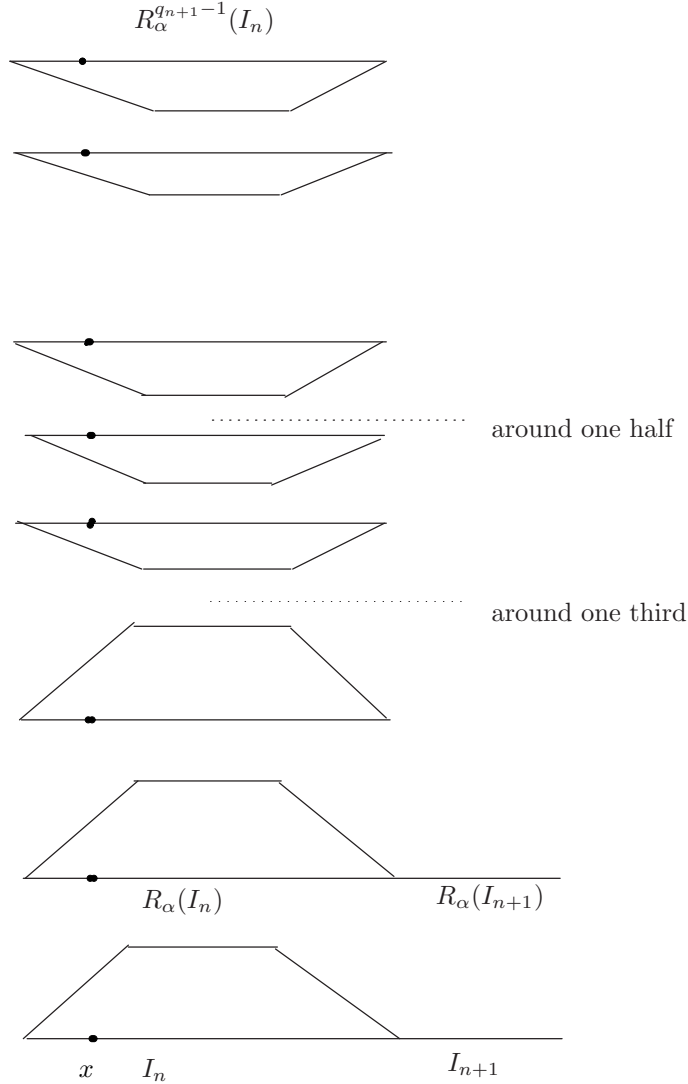

For $n$ odd, we put
$$
T_n=-\sum_{i=0}^{j_{n+1}-1}\chi_{R_\alpha^i(I_n)}+2^{-1}
\sum_{i=j_{n+1}'}^{q_{n+1}-1}\chi_{R_\alpha^i(I_n)}.$$
In any case, as is indicated in the figure, we have
\begin{equation}
 \label{ee}
T_n^{(q_{n+1})}(x)=\sum_{i=0}^{q_{n+1}-1}T_n(R_\alpha^i(x))=0 \ 
\mbox{ for any }\ x\in I_n.
\end{equation}

Let us show how the Rochlin tower for the interval $I_{n+1}\cup I_{n+2}$
is obtained from the previous one for $I_n\cup I_{n+1}$. See Figure 5
for $n$ even.
On the right side of the figure, the positive orbits of $I_{n+1}$ 
form a small tower over $I_{n+1}$ of hight $q_n$. 
Each  is mapped by $R_\alpha$ to one floor upward.
The ceiling of the tower is mapped
to the leftmost small interval on the ground level. 
Again its orbit forms a tower, this time of hight $q_{n+1}$.
Its ceiling is mapped to the second left interval on the ground level,
and so forth.

To construct the Rochlin tower for $I_{n+1}\cup I_{n+2}$,
pile up all the intervals of the size
$\abs{I_{n+1}}$ in the figure
 over the tower on $I_{n+1}$ according to the dynamical
order.  We shall get a much taller tower over
$I_{n+1}$. The narrow tower over $I_{n+2}$ in the figure
is left untouched. The
resultant is the new Rochlin tower.  Thus
the function $T_{n+1}$ looks like Figure 7.

\begin{figure}
 \caption{There are $q_{n+2}=aq_{n+1}+q_n$ intervals of size $\vert I_{n+1}\vert$.}
{\unitlength 0.1in%
}%

\end{figure}
Since $a\geq 10^{10}$, the left rectangle in Figure 2 occupies almost all portion of
the circle $S^1$.  In fact the total length of
the intervals contained in the left 
rectangle is $\vert I_n\vert\cdot q_{n+1}$, while the right 
$\vert I_{n+1}\vert\cdot q_n$. We have $q_{n+1}=aq_n+q_{n-1}>10^{10}q_n$ 
and $\vert I_n\vert>10^{10}\vert I_{n+1}\vert$ (Notice the number $a$ in Figure 5). This shows
in particular
\begin{equation}
\label{e8}
\vert I_n\vert=\vert q_n\alpha\vert \nq q_{n+1}^{-1}.
\end{equation}

\begin{proposition}
 \label{p1}
We have $\displaystyle \sum_{n=1}^\infty\vert T_n\vert_\infty<\infty$.
The series $\sum_{n=1}^\infty T_n$ converges uniformly to a continuous
function $T$.
\end{proposition}

\bd By construction, $\vert T_n\vert_\infty\leq\abs{I_n}=\vert q_n\alpha\vert\nq
q_{n+1}^{-1}$. On the other hand, by (\ref{e6}), $q_n\nq Ac^n$.
\qed

\bigskip
To show that $T$ satisfies the required property, we make use of Lemma
\ref{l1}.  Given any sufficiently small $r$,
say $0<r<10^{-100}$, we only need to compare the value
$T^{(i)}(0)$ with $T^{(i)}(x)$ for a suitably chosen $x$ such that
$\vert x\vert=r$.
Let $J_n$ be a subinterval of $I_n$ bounded by $2^{-1}q_n\alpha$ and
$-2^{-1}q_{n+1}\alpha$. See Figure 6.

For $\vert2^{-1}q_{n+1}\alpha\vert\leq r\leq\vert2^{-1}q_n\alpha\vert$, we choose
$x$ such that $\vert x\vert=r$ from the interval $J_n$, and compare the value $T^{(i)}(0)$ with
$T^{(i)}(x)$ for $i=\lfloor 2^{-1}q_{n+1}\rfloor$.
\begin{figure}
\caption{}
{\unitlength 0.1in%
\begin{picture}( 35.7000,  9.4000)( 11.8000,-22.3000)%
%
\special{pn 8}%
\special{pa 1220 1900}%
\special{pa 4750 1900}%
\special{fp}%
%
\special{pn 8}%
\special{pa 3870 1760}%
\special{pa 3870 2060}%
\special{fp}%
%
\special{pn 8}%
\special{pa 4750 1760}%
\special{pa 4750 2070}%
\special{fp}%
%
\special{pn 8}%
\special{pa 1220 1730}%
\special{pa 1220 2090}%
\special{fp}%
%
\special{pn 8}%
\special{pa 2590 1900}%
\special{pa 2590 1500}%
\special{fp}%
%
\special{pn 8}%
\special{pa 3470 1900}%
\special{pa 3470 1520}%
\special{fp}%
\put(24.5000,-22.6000){\makebox(0,0)[lb]{$I_n$}}%
\put(41.7000,-22.4000){\makebox(0,0)[lb]{$I_{n+1}$}}%
\put(29.1000,-16.2000){\makebox(0,0)[lb]{$J_n$}}%
\put(38.0000,-16.8000){\makebox(0,0)[lb]{$0$}}%
\put(46.1000,-16.5000){\makebox(0,0)[lb]{$q_{n+1}\alpha$}}%
\put(11.8000,-16.6000){\makebox(0,0)[lb]{$q_n\alpha$}}%
%
\special{pn 4}%
\special{sh 1}%
\special{ar 3090 1900 16 16 0  6.28318530717959E+0000}%
\special{sh 1}%
\special{ar 3090 1900 16 16 0  6.28318530717959E+0000}%
\put(30.5000,-21.1000){\makebox(0,0)[lb]{$x$}}%
\end{picture}}%

\end{figure}
To do this, we divide $T$ into five terms
\begin{equation}\label{e9}
T=\sum_{\nu=1}^{n-2}T_\nu+T_{n-1}+T_n+T_{n+1}+\sum_{\nu=n+2}^\infty T_\nu,
\end{equation}
and for each of these five terms, say $S$, we estimate the value
of $\displaystyle S^{(i)}(x)-S^{(i)}(0)$. Thus in the rest,
we assume the following.

\begin{assumption}
 \label{a1} $x\in J_n$ and $i=\lfloor 2^{-1}q_{n+1}\rfloor$.
\end{assumption}

First of all, let us study the middle term of (\ref{e9}).

\begin{proposition}
 \label{p2}
The sign of $T_n^{(i)}(x)-T_n^{(i)}(0)$ alternates according to
 $n$, and we have
$$
\abs{T_n^{(i)}(x)-T_n^{(i)}(0)}>100^{-1}c^{-1},$$
where $c$ is the constant in (2.2).
\end{proposition}

\bd To fix the idea, assume $n$ is even. See Figure 4. 
We are summing up the value of $T$ along one half of the vertical orbit 
which starts at the bottom line. Clearly it sums up to 0 for the initial
value 0:
$T_n^{(i)}(0)=0$.  
For the initial value $x\in J_n$,  the first
terms up to hight one third are positive, while the rest
nonpositive. We are summing up the value up to one half 
the hight of the tower of Figure 4, since $i=\lfloor 2^{-1}q_{n+1}\rfloor$.
Now $x\in J_n$ implies $\vert x\vert\geq2^{-1}\vert q_{n+1}\alpha\vert$.
The value of each first one third term is
the same and bigger than or equal to $2^{-1}\vert q_{n+1}\alpha\vert$. 
On the other hand, the rest
terms are either zero or minus half of this value.
Moreover these opposing terms are smaller in number since we are summing up until
one half the hight.
Therefore we have the following very safe estimate:
$$T_n^{(i)}(x)\geq 50^{-1}q_{n+1}\vert q_{n+1}\alpha\vert\nq
50^{-1}q_{n+1}q_{n+2}^{-1}\nq 50^{-1}c^{-1}.
$$
\qed

\begin{proposition}
 \label{p3}
We have\ \ $T_{n-1}^{(i)}(x)-T_{n-1}^{(i)}(0)=0$.
\end{proposition}

\bd We are going to show that if $x\in J_n$ and 
$i=\lfloor 2^{-1}q_{n+1}\rfloor$,
then $T_{n-1}^{(i)}(x)-T_{n-1}^{(i)}(0)=0$. In order to utilize
the previous figures, we shift the number by one. So we assume
$x\in J_{n+1}$ and $i=\lfloor 2^{-1}q_{n+2}\rfloor$, and show 
$T_{n}^{(i)}(x)-T_{n}^{(i)}(0)=0$. Thus $x$ as well as $0$ lies in
$I_{n+1}$ in Figures 5 and 6. We shall compare their orbits up to
$\lfloor2^{-1}q_{n+2}\rfloor$, half of the number of the
intervals of size $\vert I_{n+1}\vert$ in Figure 5.
 Their first $q_n$ orbits are above $I_{n+1}$ on that figure.
The values of $T_n$ sum up to 0, since $T_n$ vanishes there.
Then they come to the leftmost small interval in the bottom line. They climb the
tower and from the top falls down to the 2nd leftmost interval. At this moment, the
values of $T_n$ of both orbits sum up to 0, by (\ref{ee}).
They repeat these processes until they come to the middle of $I_n$.
At the last stage, both orbits climb up from some points in the 
middle part and stop at certain (same) hight.
But there the function $T_n$ is flat (Figure 4). Therefore the sums
for $0$ and $x$ are exactly the same. \qed

\begin{proposition}
 \label{p4} Let $x$ and $i$ be as in Assumption \ref{a1}.
The number $T_{n+1}^{(i)}(x)-T_{n+1}^{(i)}(0)$ has the same sign as
$T_n^{(i)}(x)-T_n^{(i)}(0)$.
\end{proposition}
\begin{figure}
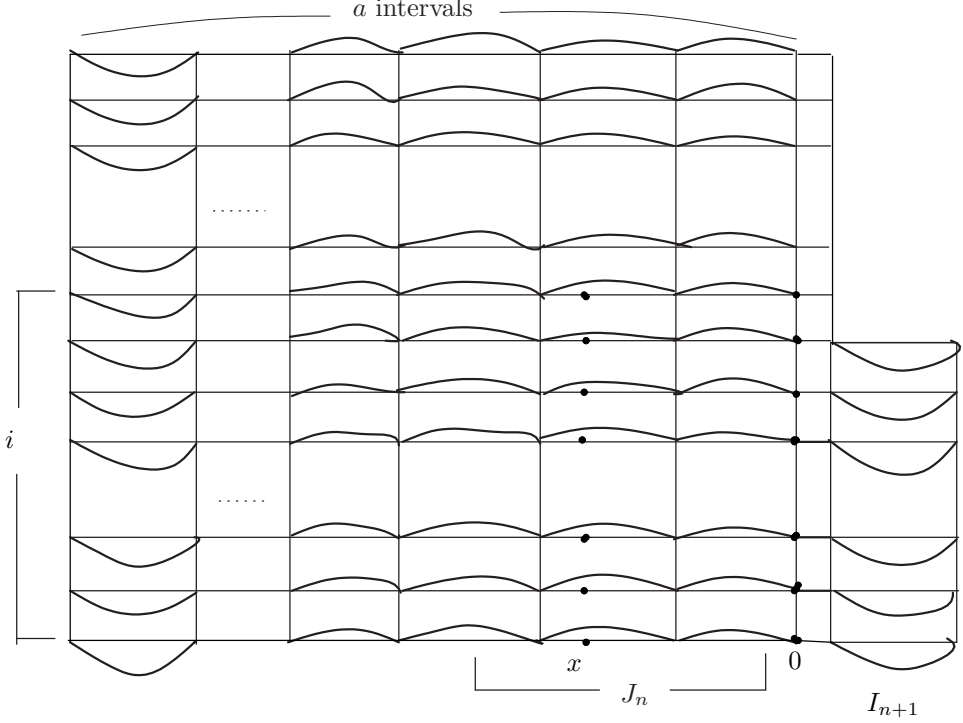

 \caption{$T_{n+1}^{(i)}(x)$ and $T_{n+1}^{(i)}(0)$ are the sum of the function at the dotted
 points.}
{\unitlength 0.1in%
}%

\end{figure}
\bd To fix the idea, assume $n$ is even. Then by the construction of $T_n$, we have
$$
T_n^{(i)}(x)-T_n^{(i)}(0)=T_n^{(i)}(x)>0.$$
Now the graph of $T_{n+1}$ is indicated in Figure 7. 
Since $n+1$ is odd, it takes negative value on $I_{n+1},
R_\alpha(I_{n+1}),\ldots$,
until at around one third of the way, it changes the sign, to positive.
It is clear from the figure that
$T_{n+1}^{(i)}(0)=0$ and $T_{n+1}^{(i)}(x)>0$. \qed

\begin{corollary}
 \label{c1} For $S=T_{n-1}+T_n+T_{n+1}$, we have $\abs{S^{(i)}(x)-S^{(i)}(0)}>100^{-1}c^{-1}$.
\end{corollary}

Now let us consider the remaining two terms in (2.6).

\begin{proposition}
 \label{p6} Assume $x$ and $i$ be as in Assumption \ref{a1}.
For $S=\sum_{\nu=n+2}^\infty T_\nu$, we have $\abs{S^{(i)}(x)-S^{(i)}(0)}<4c^{-2}$.
\end{proposition}

\bd Recall by (\ref{e6}) and (\ref{e8}) that 
$$\vert T_\nu\vert_\infty\leq \vert q_\nu\alpha\vert\nq
q_{\nu+1}^{-1}\nq A^{-1}c^{-\nu-1}.$$
This shows that 
$$\vert S\vert_\infty\leq 2A^{-1}c^{-n-3}.$$
Since $i=\lfloor 2^{-1}q_{n+1}\rfloor\nq 2^{-1}Ac^{n+1}$,
we have $\vert S^{(i)}\vert_\infty\leq 2c^{-2}$, showing the proposition. \qed

\begin{proposition}
 \label{p7} Assume $x$ and $i$ be as in Assumption \ref{a1}.
For $S=\sum_{\nu=1}^{n-2} T_\nu$, we have $\abs{S^{(i)}(x)-S^{(i)}(0)}<5c^{-2}$.
\end{proposition}

\bd For each $\nu\leq n-2$, the points $0,x\in J_n$ lie on the interval
$I_\nu\cup I_{\nu+1}$. For a point $y$ of $I_\nu\cup I_{\nu+1}$, the sum
$T_\nu^{(k)}(y)=0$ whenever $R_\alpha^k(y)$ is contained in $I_\nu\cup
I_{\nu+1}$, by (\ref{ee}).
Let ${\mathcal J}$ be the interval $I_\nu\cup I_{\nu+1}$ with the
$2^{-1}\vert q_n\alpha\vert$-neighbourhoods of the two boundary points
removed. Then if $R_\alpha^k(0)$ is contained in $\mathcal J$, 
$R_\alpha^k(x)$ is contained in $I_\nu\cup I_{\nu+1}$, since
$$
\vert R_\alpha^k(x)-R_\alpha^k(0)\vert=\vert x\vert\leq2^{-1}\vert q_n\alpha\vert.$$ 
In that case, we have $T_\nu^{(k)}(x)=T_\nu^{(k)}(0)=0$.

Now it is easy to show that the first return time of $R_\alpha$ for
$\mathcal J$ is at most $2q_{\nu+1}$.
Let $k$ be the largest integer in $\{1,2,\cdots,i\}$ such that
$R_\alpha^k(0)\in{\mathcal J}$, and let $l=i-k$. Then
$$
T_\nu^{(i)}(x)-T_\nu^{(i)}(0)
=(T_\nu^{(k)}(x)-T_\nu^{(k)}(0))+(T_\nu^{(l)}(R_\alpha^k(x))-T_\nu^{(l)}(R_\alpha^k(0)))
$$$$=T_\nu^{(l)}(R_\alpha^k(x))-T_\nu^{(l)}(R_\alpha^k(0))),$$
where $l\leq 2q_{\nu+1}$.
Since the Lipshitz constant of $T_\nu$ is 1, we have
$$
\abs{T_\nu^{(i)}(x)-T_\nu^{(i)}(0)}\leq  l\vert x\vert\leq
2q_{\nu+1}\vert q_{n+1}\alpha\vert\leq 4c^{-n+\nu}.$$
Summing up for $1\leq\nu\leq n-2$, we get Proposition \ref{p7}. \qed

\medskip
{\sc End of the proof of Theorem \ref{t2}.} For $x$ and $i$ as in
Assumption \ref{a1}, we have shown that
$$
\abs{T^{(i)}(x)-T^{(i)}(0)}>100^{-1}c^{-1}-4c^{-2}-5c^{-2}>200^{-1}c^{-1},$$
where the last inequality follows from $c>10^{10}$.
Now we have shown that the assumption of Lemma \ref{l1} is met for
$\delta=3^{-1}\min\{10^{-100},200^{-1}c^{-1}\}$.
 This shows Theorem \ref{t2}. In fact, given $\epsilon>0$.
one can choose the $\delta(\epsilon)$ in
the definition of the kinematic expansiveness as
$\delta(\epsilon)=3^{-1}\min\{\epsilon,10^{-100},200^{-1}c^{-1}\}$.

\end{document}